\documentclass[11pt]{article}
\usepackage[dvips]{graphics}
\input{amssym.def}
\input{amssym.tex}

\title{\Large\bf  Method of automorphic functions\\ for an inverse problem of antiplane elasticity}
\author{\bf {\sc Y.A. Antipov}\\ 
Department of Mathematics, Louisiana State University\\
Baton Rouge LA 70803, USA}

\date{}

\newcommand{\inte}{\mathop{\rm int}\nolimits}

\newcommand{\I}{\mathop{\rm Im}\nolimits}
\newcommand{\R}{
\mathop{\rm Re}\nolimits}
\newcommand{\const}{\mbox{const}}

\newcommand{\Md}{\partial}

\newcommand{\ov}[1]{\overline{#1}}

\newcommand{\Ga}{\alpha}
\newcommand{\Gb}{\beta}
\newcommand{\Gd}{\delta}

\newcommand{\Gf}{\phi}
\newcommand{\Gvf}{\varphi}
\newcommand{\Gg}{\gamma}
\newcommand{\Gc}{\chi}

\newcommand{\Gk}{\kappa}

\newcommand{\Gl}{\lambda}
\newcommand{\Gn}{\eta}

\newcommand{\Gt}{\theta}

\newcommand{\Gs}{\sigma}

\newcommand{\Go}{\omega}
\newcommand{\Gx}{\xi}

\newcommand{\Gz}{\zeta}

\newcommand{\GF}{\Phi}

\newcommand{\GL}{\Lambda}

\newcommand{\GO}{\Omega}

\newcommand{\CK}{{\cal K}}

\newcommand{\beq}{\begin{equation}}
\newcommand{\eeq}{\end{equation}}
\newcommand{\barr}{\begin{eqnarray}}
\newcommand{\earr}{\end{eqnarray}}
\newcommand{\beqn}{\begin{equation*}}
\newcommand{\eeqn}{\end{equation*}}
\newcommand{\barrn}{\begin{eqnarray*}}
\newcommand{\earrn}{\end{eqnarray*}}

\newcommand{\fr}{\frac}

\begin{document}
\maketitle

\begin{abstract}

A nonlinear 
 inverse problem of antiplane elasticity for a multiply connected domain is examined. It is required to determine
 the profile of $n$ uniformly stressed inclusions when the surrounding infinite body is 
subjected to antiplane uniform shear at infinity.  A method of conformal mappings of circular multiply
connected domains is employed. The conformal map is recovered by solving consequently two Riemann-Hilbert
problems for piecewise analytic symmetric automorphic functions. For domains associated with  the first
class Schottky groups a series-form  representation of a ($3n-4$)-parametric family of conformal maps
solving the problem is discovered. Numerical results for two and three uniformly stressed inclusions are reported
and discussed.

 \end{abstract}

\setcounter{equation}{0}

\section{Introduction}

Considerable interest in inverse boundary value problems for partial differential equations has developed  since the work by Riabuchinsky {\bf(\ref{ria})} 
on recovering the boundary of a domain when a function is harmonic inside 
and the function and its  normal derivative are prescribed in the boundary. Numerous applications  of inverse problems
to  filtration and hydroaerodynamics advanced the development of  a qualitative theory and various constructive techniques including those
based on the theory of boundary  value problems of the theory of analytic functions {\bf(\ref{aks1})}. The classical work by Eshelby {\bf(\ref{esh})} 
 on determination of the shapes of curvilinear  cavities and inclusions with 
prescribed properties inspired material scientists to work on inverse problems of elasticity. Eshelby found that  if the unbounded elastic body 
is uniformly loaded at infinity and the body has an elliptic or ellipsoidal inclusion with different elastic constants, 
then the  stress field is uniform inside the inclusion.  Eshelby's conjecture that in plane and antiplane cases there are no other shapes apart from ellipses 
was proved in {\bf(\ref{sen})}. Another proof based on  the method of conformal mappings was later proposed in {\bf(\ref{ru})}. 

An inverse problem of plane  elasticity for a plane uniformly
loaded at infinity and having $n$ holes was examined by Cherepanov  {\bf(\ref{che})}. In this model, the holes boundaries are  
subjected to constant normal and tangential traction, and the boundaries have to be 
determined from the condition that the tangential normal stress is constant  in all the 
contours.  Cherepanov employed the method of conformal mappings and homogeneous Schwarz problems to recover the shapes of two symmetric holes.
A circular map from the exterior of $n$-circles onto the $n$-connected elastic domain,
integral equations, and the method of least squares for their numerical solution was proposed in {\bf(\ref{vig})}.
An explicit representation in terms of the Weierstrass elliptic function for the profile
of an inclusion in the case of a doubly periodic structure was found in {\bf(\ref{gra})}.
Recently,  the theory of the Cherepanov problem for $n$ inclusions was advanced  in {\bf(\ref{ant1})} by  developing a method of the Riemann-Hilbert problem on a Riemann
surface. It was shown   that for any $n\ge 1$ there always exists a set of
the loading parameters which generate inadmissible poles of the solution.

For the antiplane inverse problem on reconstructing the boundaries of two  symmetric uniformly stressed inclusions 
the Weierstrass zeta function and the Schwarz-Christoffel formula were found  to be effective in {\bf(\ref{kan})}.
A method of Laurent series and a conformal mapping from an annulus to a doubly connected domain
 to recover the profile of two inclusions with uniform stresses was
 applied in {\bf(\ref{wan})}. Different  numerical approaches for inverse antiplane problems were employed in
  {\bf(\ref{liu})}, {\bf(\ref{dai})}.

Two methods for nonlinear inverse problems on supercavitating
flow past $n$ hydrofoils were proposed in  {\bf(\ref{ant2})}, {\bf(\ref{ant3})}. Both methods are based on the existence theorem  {\bf(\ref{kel})}, {\bf(\ref{cou})} of a conformal map of 
an $n$ connected parametric slit or circular domain into the $n$-connected physical domain.
These methods express the conformal map in terms of the solutions of two Riemann-Hilbert problems
in a multiply connected  canonical domain. In the  first method, the Riemann-Hilbert problems are set in $n$ slits
and reduce to two Riemann-Hilbert problems on a genus-$n$ Riemann surface. The parametric domain for the second technique
is the exterior of $n$ circles. The method employs  linear fractional transformations,  a symmetry transformation, and the Schottky groups {\bf(\ref{for})} and leads to
two Riemann-Hilbert problems of the theory of symmetric automorphic functions {\bf(\ref{chi})}, {\bf(\ref{sil})}, {\bf(\ref{ant4})} {\bf(\ref{ant3})}, {\bf(\ref{ant5})}.
Alternatively, the problems could be solved by using the theory of the Hilbert problems in multiply connected circular domains   
{\bf(\ref{aks2})}, {\bf(\ref{ale})}, {\bf(\ref{mit})}, {\bf(\ref{kaz})}. We also note that numerous  conformal maps of $n$-connected canonical
domains into physical domains exists in the literature when the boundary conditions of the problem allows for bypassing the Hilbert
problem in the circular domain. Examples include a map from a multiply connected circular domain into a multiply connected polygonal region 
 in a series form  {\bf(\ref{del})} and in terms of the Schottky-Klein prime function of the Schottky group associated with the circular domain  {\bf(\ref{cro})}.

In this paper we aim to propose an exact method of conformal mappings and the Riemann-Hilbert problem
of the theory of automorphic functions for the inverse antiplane problem on  $n$ 
 inclusions. The inclusions may have different shear moduli and are in ideal contact with the surrounding elastic matrix 
 subjected at infinity to uniform antiplane shear $\tau_{13}=\tau_1^\infty$
and  $\tau_{23}=\tau_2^\infty$. The profiles of the inclusions are not prescribed and have to be determined
from the condition that the stress field inside all the inclusions is uniform, $\tau_{13}=\tau_1$
and  $\tau_{23}=\tau_2$. 

In Section 2, we formulate the problem, map the exterior of $n$ circles into the exterior of $n$ uniformly stressed inclusions and 
reduce the problem of  determination of the conformal map to two inhomogeneous Schwarz problems
to be solved consecutively.  In Section 3, we convert the Schwarz problems to  two Riemann-Hilbert problems for piecewise analytic
symmetric automorphic functions. For their solution we employ a qusiautomorphic analogue of the Cauchy kernel  {\bf(\ref{chi})},  {\bf(\ref{ant4})}. Note
that the Cauchy kernel analogue  {\bf(\ref{ant6})} in terms of the Schottky-Klein prime function of the Schottky group could also be employed. 
In Section 4, for the first class Schottky groups {\bf(\ref{bur})}, we write down a series representation of a family of conformal maps solving the problem. The family has 
$3n-4$ free parameters and should satisfy the natural restriction that the inclusions contours cannot overlap.   
We also give some sample profiles of two and three uniformly stressed inclusions.
In Appendix, for completeness, we examine the case $n=1$ and show that the profile of a single uniformly stressed inclusion is an ellipse.

\vspace{.1in}

\setcounter{equation}{0}

\section{Formulation}\label{s2}

Consider the following problem of antiplane elasticity.

{\sl 
Let  $D_0$, $D_1$, $\ldots,D_{n-1}$ be $n$ finite inclusions in an infinite isotropic solid. The shear moduli of the inclusions and
the solid $D^e={\Bbb R}^2\setminus D$  $(D=\cup_{j=0}^{n-1}D_j)$ are taken to be $\mu_j$ and $\mu$, respectively.
It is assumed that the inclusions are in ideal contact  with 
the matrix, and the whole solid $D^e\cup D$ is in a state of antiplane shear  due 
to constant shear stresses applied at infinity, $\tau_{13}=\tau_1^\infty$, $\tau_{23}=\tau_2^\infty$. 
It is  aimed to determine the boundaries of the inclusions, $L_j$, such that the stresses
$\tau_{13}$ and $\tau_{23}$ are constant in all the inclusions $D_j$, $\tau_{13}=\tau_1$, $\tau_{23}=\tau_2$,
$j=0,1,\ldots,n-1$. }

\vspace{2mm}

Let $u$ and $u_j$ be the $x_3$-components of the displacement vectors for the body $D^e$ and the inclusions $D_j$, respectiely. 
Then $\tau_{k3}=\mu\Md u/\Md x_k$
($k=1,2$), $(x_1,x_2)\in D^e$, and $\tau_{k3}=\mu_j\Md u/\Md x_k$
($k=1,2$), $(x_1,x_2)\in D_j$, $j=0,1,\ldots,n-1$.  At infinity, the displacement $u$ is growing as 
\beq
u\sim \mu^{-1}(\tau_1^\infty x_1+\tau_2^\infty x_2)+\const,\quad  x_1^2+x_2^2\to\infty.
\label{2.1}
\eeq
Due to the fact that the stresses $\tau_{12}$ and $\tau_{13}$ are constant in the inclusions, the  $x_3$-displacements
$u_j$ for ($x_1,x_2)\in D_j$ are linear functions
\beq
u_j=\mu_j^{-1}(\tau_1 x_1+\tau_2 x_2)+d_j', \quad  (x_1,x_2)\in D_j, \quad j=0,1,\ldots,n-1,
\label{2.2}
\eeq
and $d_j'$ are real constants.

Let $v$ and $v_j$ be the harmonic conjugates of the  harmonic functions   $u$ and $u_j$ in the domains $D^e$ and $D_j$, respectively. 
Denote $z=x_1+ix_2$.
Then
$\Gf(z)=u(x_1,x_2)+iv(x_1,x_2)$ and $\Gf_j(z)=u_j(x_1,x_2)+iv_j(x_1,x_2)$ are analytic  functions in the corresponding domains $D^e$ and $D_j$.
The boundary conditions of ideal contact imply
that the traction component $\tau_{\nu3}$ and  the $x_3$-component of the displacement 
are continuous  through the contours $L_j$,
\beq
\mu\fr{\Md u}{\Md \nu}=\mu_j\fr{\Md u_j}{\Md \nu}, \quad u=u_j, \quad (x_1,x_2)\in L_j, \quad j=0,1,\ldots,n-1,
\label{2.3}
\eeq
where $\fr{\Md}{\Md \nu}$ is the normal derivative. In terms of the functions $\Gf_j(z)$ the ideal contact boundary conditions can be written as
\beq
\fr{\Gk_j+1}{2}\Gf_j(z)-\fr{\Gk_j-1}{2}\ov{\Gf_j(z)}=\Gf(z)+ib_j, \quad z\in L_j, \quad j=0,1,\ldots,n-1,
\label{2.4}
\eeq
where $\Gk_j=\mu_j/\mu$, and $b_j$ are real constants. The equivalence of the boundary conditions (\ref{2.3})
 and (\ref{2.4}) follows from the Cauchy-Riemann condition $\fr{\Md u}{\Md\nu}=\fr{\Md v}{\Md s}$,
where $\fr{\Md}{\Md s}$ is the tangential derivative. Since the functions $u_j$ are known everywhere in the domains $D_j$,
the functions  $\Gf_j(z)$ are defined up to arbitrary constants and given by by
\beq
\Gf_j(z)=\fr{\bar\tau z}{\mu_j}+d_j, \quad z\in D_j, \quad j=0,1,\ldots,n-1,
\label{2.5}
\eeq
where $\bar\tau=\tau_1-i\tau_2$, $d_j=d_j'+id_j''$, $d_j''$ are real constants.
In view of the relations  (\ref{2.5}), instead of the function $\Gf(z)$, it is convenient to deal with the function 
$f(z)=\Gf(z)-\bar\tau z/\mu$, $z\in D^e$.
The new function $f(z)$ is analytic in the domain $D^e$, satisfies the boundary condition
\beq
f(z)=\fr{1}{\Gl_j}\R\left(\fr{\bar\tau}{\mu} z\right)+d_j'+ia_j, \quad z\in L_j, \quad j=0,1,\ldots,n-1,
\label{2.6}
\eeq
and the condition at infinity
\beq
f(z)\sim \fr{(\bar\tau^\infty-\bar\tau)z}{\mu}+\const, \quad z\to\infty.
\label{2.7}
\eeq
Here, $\Gl_j=\Gk_j/(1-\Gk_j)$, $a_j=\Gk_j d_j''-b_j$ are real constants, and $\bar\tau^\infty=\tau^\infty_1-i\tau^\infty_2$. The condition (\ref{2.7}) 
is due to the relation $\Gf(z)\sim\bar\tau^\infty z/\mu +\const$, $z\to\infty$.

Let $z=\Go(\Gz)$ be a conformal map that transforms the exterior of $n$ circles $\frak L_j$ ($j=0,1,\ldots,n-1$)
into the physical domain $D^e$. Denote this $n$-connected circular domain by $\frak D^e$.  By scaling and rotation, it is always possible to achieve 
{\bf(\ref{kel})}, {\bf(\ref{cou})} that one of
the circles say, $\frak L_0$, is of unit radius and centered at the origin and, in addition, the center of another
circle say, $\frak L_1$, falls in the real axis. Let the circular map meet the condition $\Go(\infty)=\infty$. In this case, if the original problem
has a unique solution, then the radius of $\frak L_1$, the complex centers and the radii 
of the rest $n-2$ circles cannot be selected arbitrarily. 

In the vicinity of the infinite point the conformal map can be represented as
\beq
\Go(\Gz)=c_{-1}\Gz+c_0+\sum_{j=1}^\infty \fr{c_j}{\Gz^j},
\label{2.8}
\eeq
where  $c_{-1}=c_{-1}'+ic_{-1}''$.
Denote $f(\Go(\Gz))=F(\Gz)$. From the boundary condition (\ref{2.6}) we deduce  that the functions $F(\Gz)$
and $\Go(\Gz)$ satisfy the following two Schwarz problems to be solved consecutively.

{\it Find two functions $F(\Gz)$ and $\Go(\Gz)$ analytic in the domain $\frak D^e$ and continuous up to the boundary
$\frak L=\cup_{j=0}^{n-1} \frak L_j$ such that
\beq
\I F(\Gz)=a_j, \quad \Gz\in \frak L_j, \quad j=0,1\ldots,n-1,
\label{2.10}
\eeq
and
\beq
\R \left[\fr{\bar\tau}{\mu}\Go(\Gz)\right]=\Gl_j[\R F(\Gz)-d_j'], \quad \Gz\in \frak L_j, \quad j=0,1\ldots,n-1.
\label{2.11}
\eeq
At the infinite point, both of the  functions have a simple pole,
\beq
F(\Gz)\sim \fr{\bar\tau^\infty-\bar\tau}{\mu}c_{-1}\Gz, \quad \Go(\Gz)\sim c_{-1}\Gz, \quad \Gz\to\infty.
\label{2.13}
\eeq
In addition, the function $\Go: \frak L_j\to L_j$ ($j=0,\ldots, n-1$) has to be univalent, and the interiors of the images of the contours $\frak L_j$, the domains $D_j$, are disjoint sets.
}

\setcounter{equation}{0}

\section{Riemann-Hilbert problems of the theory of automorphic functions}

\subsection{Setting}

Let $\frak G$ be the symmetry group of the circular line $\frak L=\frak L_0\cup\ldots\cup\frak L_{n-1}$ generated by
the linear transformations $\Gs_j=T_j T_0(\Gz)$, $j=0,1,\ldots,n-1$, where $T_j$ are linear fractional transformations
\beq
T_j(\Gz)=\Gz_j+\fr{r_j^2}{\bar\Gz-\bar\Gz_j}, \quad j=0,1,\ldots,n-1.
\label{3.1}
\eeq
Here, $r_j$ and $\Gz_j$ are the radius and the center of the circle $\frak L_j$.
The transformation $\Gs_j(\Gz)$  ($j=1,\ldots,n-1$) maps the exterior domain $\frak D^e$
into the exterior of $\Gs_j(\frak L_m)\subset \inte {\frak L_j}$, $j,m=1,\ldots,n-1$. Denote by $\tilde{\frak D}^e=T_0(\frak D^e)$ 
the exterior of the circles $T_0(\frak L_j)$ inside $\frak L_0$.
The domains $\frak D^e$  and $\tilde{\frak D}^e$ and the contour $\frak L$ comprises a fundamental region $\frak F_{\frak G}$ 
of the group $\frak G$: $\frak F_{\frak G}=\frak D\cup \tilde\frak D^e\cup\frak L$. The group $\frak G$ is a symmetry
Schottky group {\bf(\ref{for})}. It consists of the identical map $\Gs_0(\Gz)=T_0T_0(\Gz)=\Gz$ and all possible compositions 
of the generators $\Gs_j$ and the inverse maps $\Gs_j^{-1}=T_0T_j$ $j=1,2,\ldots,n-1$. Therefore, each element
of the group $\frak G$ is a composition of an even number of the symmetry maps $T_j(\Gz)$ ($j=0,1,\ldots,n-1$),
$$
\Gs=T_{k_1}T_{k_2}\ldots T_{k_{2s-1}}T_{k_{2s}}, 
$$
\beq
 k_2\ne k_1, k_3\ne k_2,\ldots, k_{2s}\ne k_{2s-1}, \quad k_1,k_2,\ldots k_{2s}=0,1,\ldots,n-1.
\label{3.2}
\eeq
The region $\GO=\cup_{\Gs\in\frak G}\Gs(\frak F_{\frak G})$ is invariant with respect to the group $\frak G$:
$\Gs(\GO)=\GO$ for any $\Gs\in\frak G$, where $\GO=\bar{\Bbb C}\setminus\GL$, $\bar{\Bbb C}=\Bbb C\cup\{\infty\}$,
and $\GL$ is the set of  all limit points of the group $\frak G$.  If $n=2$, then the set $\GL$ consists of two points, while
for $n\ge 3$, the number of limit points is infinite. All maps of the group $\frak G$ are 
linear fractional transformations
\beq
\Gs(\Gz)=\fr{\Ga_\Gs\Gz+\Gb_\Gs}{\Gg_\Gs\Gz+\Gd_\Gs}, \quad \Ga_\Gs\Gd_\Gs-\Gg_\Gs\Gb_\Gs\ne 0, 
\label{3.3}
\eeq
and $ \Gg_\Gs\ne 0$ if $\Gs\ne\Gs_0$. It is assumed that
the series
\beq
\sum_{\Gs\in{\frak G}\setminus\Gs_0}\fr{|\Ga_\Gs\Gd_\Gs-\Gb_\Gs\Gg_\Gs|}{|\Gg_\Gs|^2}
\label{3.4}
\eeq
is convergent.
If $\frak G$ is a first class group {\bf(\ref{bur})}, that is if the condition (\ref{3.4}) is satisfied, then
each $\frak G$-automorphic function is representable as a series with simple fractions as its elements. 
By the sufficient Schottky condition  {\bf(\ref{sch})} $\frak G$ is a first class group if the
domain $\frak D$ can be split into a union of triply or doubly connected domains by circles which do not
intersect each other and the circles $\frak L_j$ ($ j = 0, 1, \ldots, n-1$). 
Remark that $\frak G$ is a first class Schottky group if it is associated with a doubly or
 triply connected domains, or a  circular multiply connected domain when the centers of the
circles $\frak L_j$  lie in a straight line, or it satisfies the Aksentiev condition {\bf(\ref{aks2})}. An example of the domain $\frak D$ that generates a
 symmetry Schottky group for which the series (3.5) is divergent and the corresponding Poincar\'e series of dimension $-2$ is not absolutely convergent is given in {\bf(\ref{ant4})}.
  In what follows, it is assumed that the  domain $\frak D$ obeys the sufficient conditions which guarantee the convergence of the series (\ref{3.4}). This justifies the change of order of summation used in the representation of a quasiautomorphic analogue of the Cauchy kernel.

Introduce next two functions, $\GF_1(\Gz)$ and   $\GF_2(\Gz)$, analytic in the domain $\frak D^e$ by
\beq
\GF_1(\Gz)=F(\Gz)-c\Gz, \quad \GF_2(\Gz)=\fr{i\bar\tau}{\mu}[\Go(\Gz)-c_{-1}\Gz], \quad \Gz\in{\frak D^e},
\label{3.5}
\eeq
where
\beq
c=\fr{\bar\tau^\infty-\bar\tau}{\mu}c_{-1}.
\label{3.6}
\eeq
To extend their definition inside the domain $\frak D$, we set
$$
\GF_m(\Gz)=\ov{\GF_m(T_0(\Gz))}, \quad \Gs\in T_0(\frak D^e),
$$
\beq
\GF_m(\Gz)=\GF_m(\Gs^{-1}(\Gz)), \quad \Gz\in\Gs(\frak D^e\cup T_0(\frak D^e)), \quad \Gs\in\frak G.
\label{3.7}
\eeq
Then $\GF_m(\Gz)$ are piecewise meromorphic and $\frak G$-automorphic functions which satisfy the symmetry condition
\beq
\ov{\GF_m(T_j(\Gz))}=\GF_m(\Gz), 
\quad \Gz\in T_j(\frak D^e)=\Gs_j(T_0(\frak D^e)), \quad j=1,2,\ldots n-1.
\label{3.8}
\eeq
All circles $\Gs(\frak L)$ including $\frak L$ are discontinuity lines for the functions $\GF_m(\Gz)$.
Let $\GF_m^+(\Gx)$ and $\GF_m^-(\Gx)$ be the boundary values of the functions $\GF_m(\Gz)$
from the interior and the exterior of the circles $\Gs(\frak L)$, $\Gs\in\frak G$, respectively. Then the functions
$\GF_m(\Gz)$ solve the following two Riemann-Hilbert problems.

{\sl Find all piecewise analytic and $\frak G$-automorphic functions bounded at infinity
which meet the symmetry condition (\ref{3.8}) and satisfy the linear relation in the circles $\frak L_j$ 
\beq
\GF^+_m(\Gx)-\GF^-_m(\Gx)=g_{mj}(\Gx), \quad \Gx\in \frak L_j, \quad  j=0,1,\ldots,n-1,\quad m=1,2,
\label{3.9}
\eeq
where
\beq
g_{1j}(\Gx)=2i[\I(c\Gx)-a_j], \quad g_{2j}(\Gx)=-2i\R\left\{\fr{\Gk_j}{1-\Gk_j}[F(\Gx)-d_j']-\fr{\bar\tau}{\mu}c_{-1}\Gx\right\}.
\label{3.10}
\eeq
 }

Due to formulas (\ref{3.5}) to determine the conformal mapping $\Go(\Gz)$, it is required to find the function $\GF_2(\Gz)$.
This can be done only if the first function $\GF_1(\Gz)$ is known. To solve the two Riemann-Hilbert problems (\ref{3.9}), we 
employ the following singular integral:
\beq
\Psi(\Gz)=\fr{1}{2\pi i}\int_{\frak L} \CK(\Gz,\Gn)g(\Gn)d\Gn,
\label{3.11}
\eeq
where $\CK(\Gz,\Gn)$ is the series {\bf(\ref{chi})},  {\bf(\ref{ant4})}
\beq
\CK(\Gz,\Gn)=\sum_{\Gs\in\frak G}\left(\fr{1}{\Gs(\Gn)-\Gz}-\fr{1}{\Gs(\Gn)-\Gz_*}\right)\Gs'(\Gn),
\label{3.12}
\eeq
$\Gz_*\in\frak D^e$ is an arbitrary fixed point, and $g(\Gn)$ is a density. 
Alternatively, because of the identity
\beq
\left(\fr{1}{\Gs(\Gn)-\Gz}-\fr{1}{\Gs(\Gn)-\Gz_*}\right)\Gs'(\Gn)=\fr{1}{\Gn-\Gs^{-1}(\Gz)}-\fr{1}{\Gn-\Gs^{-1}(\Gz_*)},
\label{3.12'}
\eeq
the kernel $\CK(\Gz,\Gn)$ can be written as
\beq
\CK(\Gz,\Gn)=\sum_{\Go\in\frak G}\left(\fr{1}{\Gn-\Go(\Gz)}-\fr{1}{\Gn-\Go(\Gz_*)}\right).
\label{3.12''}
\eeq
Here, since $\Gs\in\frak G$ implies $\Gs^{-1}\in\frak G$, we made the substitution  $\Go=\Gs^{-1}$.
The absolute convergence of the series in (\ref{3.12'}) and (\ref{3.12''}) is guaranteed by the convergence of the series (\ref{3.4}).
This circumstance allows for changing the order of summation and integration in the expression (\ref{3.11}).
Since the identity transformation $\Gs_0\in{\frak G}$, the series admits
the representation
\beq
\CK(\Gz,\Gn)=\fr{1}{\Gn-\Gz}+\CK_0(\Gz,\Gn),
\label{3.13} 
\eeq
where $\CK_0(\Gz,\Gn)$ is an analytic function of $\Gz$ in the domain $\GO$
\beq
\CK_0(\Gz,\Gn)=-\fr{1}{\Gn-\Gz_*}+\sum_{\Gs\in\frak G\setminus\Gs_0}\left(\fr{1}{\Gn-\Gs(\Gz)}-\fr{1}{\Gn-\Gs(\Gz_*)}\right).
\label{3.13'}
\eeq
For our next step, we need to use the following property of the kernel $\CK(\Gz,\Gn)$  {\bf(\ref{chi})},  {\bf(\ref{ant4})}:
\beq
\CK(\Gs(\Gz),\Gn)=\CK(\Gz,\Gn)+\Gc_\Gs(\Gn), \quad \Gs\in\frak G,
\label{3.14}
\eeq
where
\beq
\Gc_\Gs(\Gn)=\CK(\Gs(\Gz_*),\Gn).
\label{3.15}
\eeq
The two relations (\ref{3.13}) and (\ref{3.14}) classify the function $\CK(\Gz,\Gn)$
as a quasiauthomorphic  analogue of the Cauchy kernel.

\subsection{Solution $\GF_1(\Gz)$ of the first Riemann-Hilbert problem}

We next claim that under a certain choice of the real constants $a_j$ the  function
\beq
\GF_1(\Gz)=\Psi_1(\Gz)+\ov{\Psi_1(T_0(\Gz))}+C_1,
\label{3.16}
\eeq
provides the general solution of the first Riemann-Hilbert problem (\ref{3.9}). Here, $C_1$ is an arbitrary real constant and
\beq
\Psi_1(\Gz)=\fr{1}{2\pi}\sum_{j=0}^{n-1}\int_{\frak L_j}[\I(c\Gn)-a_j]\CK(\Gz,\Gn)d\Gn.
\label{3.17}
\eeq
Apparently, the function $\GF_1(\Gz)$ satisfies the symmetry condition $\GF_1(\Gz)=\ov{\GF_1(T_0(\Gz))}$, $\Gs\in T_0(\frak D^e)$.
Because of the quasiautomorphicity of the kernel $\CK(\Gz,\Gn)$, in general, the function $\GF_1(\Gz)$ is not automorphic. 
However, by a certain choice of the constants $a_j$, it is possible to satisfy the condition $\GF_1(\Gz)=\GF_1(\Gs_j(\Gz))$, $j=1,2,\ldots,n-1$, and 
make the solution automorphic. Indeed, the relation (\ref{3.14}) implies
\beq
\CK(\Gs_j(\Gz),\Gn)=\CK(\Gz,\Gn)+\Gc_{\Gs_j}(\Gn), \quad \Gc_{\Gs_j}(\Gn)=\CK(\Gs_j(\Gz_*),\Gn), \quad j=1,2,\ldots,n-1.
\label{3.18}
\eeq
Hence 
\beq
\Psi_1(\Gs_j(\Gz))=\Psi_1(\Gz)+e_j, \quad j=1,2,\ldots,n-1,
\label{3.19}
\eeq
where
\beq
e_j=\fr{1}{2\pi}\sum_{l=0}^{n-1}\int_{\frak L_l}[\I(c\Gn)-a_l]\Gc_{\Gs_j}(\Gn)d\Gn.
\label{3.20}
\eeq
On writing the relation (\ref{3.19}) as $\Psi_1(\Gz)=\Psi_1(\Gs^{-1}_j(\Gz))+e_j$,  $j=1,2,\ldots,n-1,$ there is no difficulty 
in verifying that
\beq
\ov{\Psi_1(T_0\Gs_j(\Gz))}=\ov{\Psi_1(\Gs^{-1}_jT_0(\Gz))}=\ov{\Psi_1(T_0(\Gz))}-\ov{e_j}.
\label{3.21}
\eeq
In view of (\ref{3.19}) and (\ref{3.21}) we obtain
\beq
\GF_1(\Gs_j(\Gz))=\GF_1(\Gz)+e_j-\ov{e_j}, \quad j=1,2,\ldots,n-1.
 \label{3.22}
 \eeq
 Hence the solution $\GF_1(\Gz)$ is a $\frak G$-automorphic function if and only if $\I e_j=0$, $j=1,2,\ldots,n-1$, 
 that is
 \beq
 \I\left\{\sum_{l=0}^{n-1}\int_{\frak L_l}[a_l-\I(c\Gn)]\Gc_{\Gs_j}(\Gn)d\Gn\right\}=0, \quad j=1,2,\ldots,n-1.
 \label{3.23}
 \eeq
The integrals 
\beq
\int_{\frak L_l}\Gc_{\Gs_j}(\Gn)d\Gn =\sum_{\Gs\in\frak G}\int_{\frak L_l}\left(\fr{1}{\Gn-\Gs\Gs_j(\Gz_*)}-\fr{1}{\Gn-\Gs(\Gz_*)}\right)d\Gn, 
\quad j=1,2,\ldots,n-1,
\label{3.24}
\eeq
can be evaluated by the theory of residues. Assume first that $l=0$. If $\Gs=\Gs_0$, then $\Gs\Gs_j(\Gz_*)\in\inte \frak L_j$ ($j=1,2,\ldots, n-1$)
and
$\Gs(\Gz_*)=\Gz_*\in\frak D^e$.  Hence,
\beq
\int_{\frak L_0}\left(\fr{1}{\Gn-\Gs\Gs_j(\Gz_*)}-\fr{1}{\Gn-\Gs(\Gz_*)}\right)d\Gn=0.
\label{3.25}
\eeq
If $\Gs=\Gs_j^{-1}$ and since $\Gs_j^{-1}=T_0T_j$, then  $\Gs\Gs_j(\Gz_*)\in\frak D^e$ and $\Gs(\Gz_*)=T_0T_j(\Gz_*)\in\inte\frak L_0$. Thus,
\beq
\int_{\frak L_0}\left(\fr{1}{\Gn-\Gs\Gs_j(\Gz_*)}-\fr{1}{\Gn-\Gs(\Gz_*)}\right)d\Gn=-2\pi i.
\label{3.26}
\eeq
Let now $\Gs\ne\Gs_0$, $\Gs\ne\Gs_j^{-1}$, and $\Gs=T_k\ldots T_\nu$. If  $k\ne 0$, then
\beq
\int_{\frak L_0}\fr{d\Gn}{\Gn-\Gs\Gs_j(\Gz_*)}=\int_{\frak L_0}\fr{d\Gn}{\Gn-\Gs(\Gz_*)}=0.
\label{3.27}
\eeq
In the case $k=0$ we have $\Gs\Gs_j(\Gz_*)\in\inte\frak L_0$ and  $\Gs(\Gz_*)\in\inte\frak L_0$. That is why
\beq
\int_{\frak L_0}\fr{d\Gn}{\Gn-\Gs\Gs_j(\Gz_*)}=\int_{\frak L_0}\fr{d\Gn}{\Gn-\Gs_j(\Gz_*)}=2\pi i.
\label{3.28}
\eeq
Summing up the results obtained we evaluate the integral (\ref{3.24}) for $l=0$
\beq
\int_{\frak L_0}\Gc_{\Gs_j}(\Gn)d\Gn =-2\pi i
\quad j=1,2,\ldots,n-1.
\label{3.30}
\eeq
Assume next that $l=1,2,\ldots,n-1$ and evaluate the integrals (\ref{3.24}). If $\Gs=\Gs_0$, then $\Gs\Gs_j(\Gz_*)\in\inte\frak L_j$, while
$\Gs(\Gz_*)=\Gz_*\notin \inte\frak L_j$. This implies 
\beq
\int_{\frak L_l}\fr{d\Gn}{\Gn-\Gs\Gs_j(\Gz_*)}=\left\{
\begin{array}{cc}
2\pi i, & j=l,\\
0, & j\ne l,\\
\end{array}\right.
\quad
\int_{\frak L_l}\fr{d\Gn}{\Gn-\Gs_j(\Gz_*)}=0, \quad l=1,2,\ldots,n-1.
\label{3.31}
\eeq
In the case $\Gs=\Gs_j^{-1}$ we have $\Gs\Gs_j(\Gz_*)=\Gz_*\in\frak D^e$ and $\Gs(\Gz_*)\in\inte \frak L_0$, and therefore
\beq
\int_{\frak L_l}\fr{d\Gn}{\Gn-\Gs\Gs_j(\Gz_*)}=\int_{\frak L_l}\fr{d\Gn}{\Gn-\Gs(\Gz_*)}=0.
\label{3.32}
\eeq
Suppose $\Gs\ne\Gs_0$, $\Gs\ne\Gs_j^{-1}$, and $\Gs=T_l\ldots T_\nu$. This implies that $\Gs\Gs_j(\Gz_*))\in\inte\frak L_l$ and
$\Gs(\Gz_*)\in\inte\frak L_l$. Therefore
\beq
\int_{\frak L_l}\fr{d\Gn}{\Gn-\Gs\Gs_j(\Gz_*)}=\int_{\frak L_l}\fr{d\Gn}{\Gn-\Gs(\Gz_*)}=2\pi i.
\label{3.33}
\eeq
If $\Gs$ takes on the value $T_k\ldots T_\nu$ and $k\ne l$, then 
\beq
\int_{\frak L_l}\fr{d\Gn}{\Gn-\Gs\Gs_j(\Gz_*)}=\int_{\frak L_l}\fr{d\Gn}{\Gn-\Gs(\Gz_*)}=0.
\label{3.34}
\eeq 
Combining all these cases we discover
\beq
\int_{\frak L_l}
\Gc_{\Gs_j}(\Gn)d\Gn=
2\pi i\Gd_{lj}, \quad l,j=1,2,\ldots,n-1.
\label{3.35}
\eeq
On substituting the integrals (\ref{3.30}) and (\ref{3.35}) into equations  (\ref{3.23})
we determine all the constants $a_1, a_2, \ldots, a_{n-1}$
\beq
a_j=a_0+\fr{1}{2\pi}\I\sum_{l=0}^{n-1}\int_{\frak L_l}\I(c\Gn)\Gc_{\Gs_j}(\Gn)d\Gn, \quad j=1,2,\ldots,n-1.
\label{3.36}
\eeq
The constant $a_0$ remains to be free.
We thus proved that if the constants $a_j$ are chosen as in (\ref{3.36}), then $\GF_1(\Gz)$
is a $\frak G$-automorphic function. Show finally that it satisfies the Riemann-Hilbert boundary condition (\ref{3.9}).
On splitting $\frak G$ into $\Gs_0$ and $\frak G\setminus\Gs_0$ we represent the function $\Psi_1(\Gz)$
in the form
$$
\Psi_1(\Gz)=\fr{1}{2\pi}\sum_{l=0}^{n-1}\int_{\frak L_l}\left(\fr{1}{\Gn-\Gz}-\fr{1}{\Gn-\Gz_*}\right)[\I(c\Gn)-a_l]d\Gn
$$
\beq
+\fr{1}{2\pi}\sum_{l=0}^{n-1}
\sum_{\Gs\in\frak G\setminus\Gs_0}
\int_{\frak L_l}\left(\fr{1}{\Gn-\Gs(\Gz)}-\fr{1}{\Gn-\Gs(\Gz_*)}\right)[\I(c\Gn)-a_l]d\Gn.
\label{3.37}
\eeq
Passing to the limit $\Gz\to\Gx\in\frak L_j$ and employing the Sokhotski-Plemelj formulas we obtain
\beq
\Psi_1^-(\Gx)=-\fr{i}{2}[\I(c\Gx)-a_j]+\Psi_1(\Gx).
\label{3.38}
\eeq
In a similar fashion we next analyze the function $\ov{\Psi_1(T_0(\Gz))}$. We have 
\beq
\lim_{\Gz\to\Gx\in\frak L_j, \Gz\in\frak D^e}\ov{\Psi_1(T_0(\Gz))}=
-\fr{i}{2}[\I(c\Gx)-a_j]
+\ov{\Psi_1(T_0(\Gx))},
\label{3.38'}
\eeq
and then the limit values $\GF_1^\pm(\Gx)$of the solution become
\beq
\GF_1^\pm(\Gx)=\pm i[\I(c\Gx)-a_j]+\Psi_1(\Gx)+\ov{\Psi_1(T_0(\Gx))}, \quad \Gx\in \frak L_j.
\label{3.39}
\eeq
This verifies that the $\frak G$-automorphic symmetric function $\GF_1(\Gz)$ bounded at infinity and given by (\ref{3.16})
with the constants $a_j$ defined by (\ref{3.36}) solves the first Riemann-Hilbert problem (\ref{3.9}). Any other solution of this problem 
differs from the function (\ref{3.16}) by a constant. This may be proved in a manner standard in the theory of boundary value problems of the theory
of analytic functions {\bf(\ref{gak})}.

\subsection{The function $\GF_2(\Gz)$ }

We begin with rewriting the second Riemann-Hilbert problem (\ref{3.9}) in the following form:
\beq
\GF_2^+(\Gx)-\GF_2^-(\Gx)=-2i\tilde d_j+2ig_{2j}^\circ(\Gx), \quad \Gz\in\frak L_j, \quad j=0,1,\ldots,n-1,
\label{3.40}
\eeq
where
\beq
\tilde d_j=\fr{\Gk_j}{\Gk_j-1}d_j',\quad
g_{2j}^\circ =\R\left(\fr{\Gk_j}{\Gk_j-1}F(\Gx)+\fr{\bar\tau }{\mu}c_{-1}\Gx\right),
\label{3.41}
\eeq
and
\beq
F(\Gx)=c\Gx-i[\I(c\Gx)-a_j]+\Psi_1(\Gx)+\ov{\Psi_1(T_0(\Gx))}, \quad \Gx\in \frak L_j.
\label{3.42}
\eeq
The Riemann-Hilbert problem is solved in the same fashion as the first problem for the function $\GF_1(\Gz)$.
Its solution is given by
\beq
\GF_2(\Gz)=\Psi_2(\Gz)+\ov{\Psi_2(T_0(\Gz))}+C_2,
\label{3.43}
\eeq
where $C_2$ is an arbitrary real constant and
\beq
\Psi_2(\Gz)=\fr{1}{2\pi}\sum_{j=0}^{n-1}\int_{\frak L_j}[g_{2j}^\circ(\Gn)-\tilde d_j]\CK(\Gz,\Gn)d\Gn, \quad \Gz\notin \frak L.
\label{3.44}
\eeq
This functions is a $\frak G$-automorphic function if and only if  the constants 
$\tilde d_j$ are given by
\beq
\tilde d_j=\tilde d_0+\fr{1}{2\pi}\I\sum_{l=0}^{n-1}\int_{\frak L_l} g_{2l}^\circ(\Gn)\Gc_{\Gs_j}(\Gn)d\Gn, \quad j=1,2,\ldots,n-1.
\label{3.45}
\eeq
The constant $\tilde d_0$ may be fixed arbitrarily.  
As before, it is directly verified that the function $\GF_2(\Gx)$
satisfies the symmetry condition (\ref{3.8}).  On the circles $\frak L_j$, the limit values of the function $\GF_2(\Gx)$ are determined according to the Sokhotski-Plemelj
formulas
\beq
\GF_2^\pm(\Gx)=\pm i[g_{2j}^\circ(\Gx)-\tilde d_j]+\Psi_2(\Gx)+\ov{\Psi_2(T_0(\Gx))}, \quad \Gx\in \frak L_j.
\label{3.46}
\eeq
where $\Psi_2(\Gx)$ and $\ov{\Psi_2(T_0(\Gx))}$ are the Cauchy principal values of the integrals $\Psi_2(\Gz)$ and $\ov{\Psi_2(T_0(\Gz))}$, respectively,  with $\Gz=\Gx\in\frak L_j$.

\setcounter{equation}{0} 
 
\setcounter{equation}{0}

\section{Conformal mapping. Numerical results}

\begin{figure}[t]
\centerline{
\scalebox{0.6}{\includegraphics{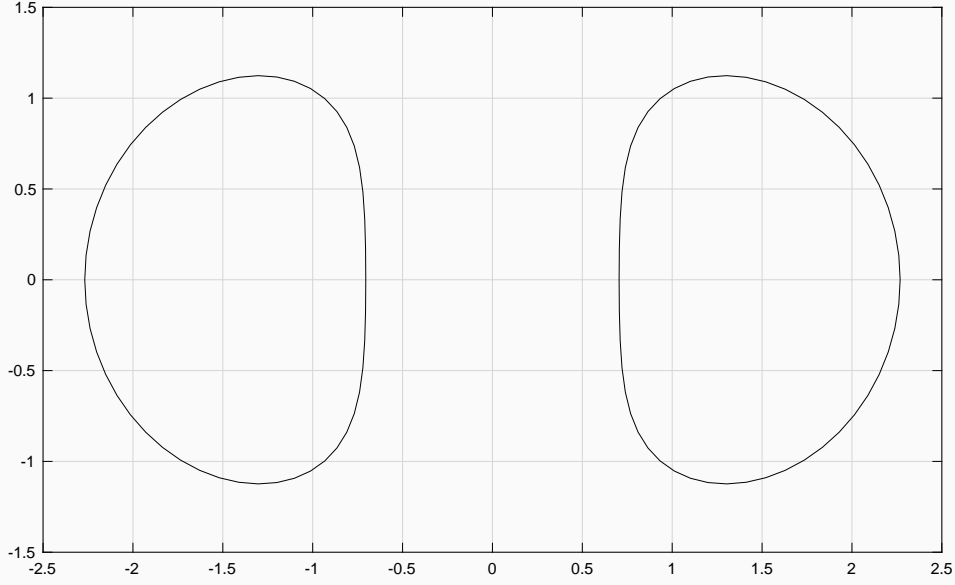}}}
\caption{Two symmetric inclusions $(n=2)$ when 
$\tau_1/\mu=2$, $\tau_1^\infty/\mu=1 $, $\tau_2=\tau_2^\infty=0$, $\Gk_0=\Gk_1=2$,   $c_{-1}=1$, $r_0=r_1=1$,  $\Gz_0=-1.5$,  $\Gz_1=1.5$,  $\Gz_*=0$.
$a_0=-a_1$, and $\tilde d_0=-\tilde d_1$.
} 
\label{fig1}
\end{figure}

\begin{figure}[t]
\centerline{
\scalebox{0.6}{\includegraphics{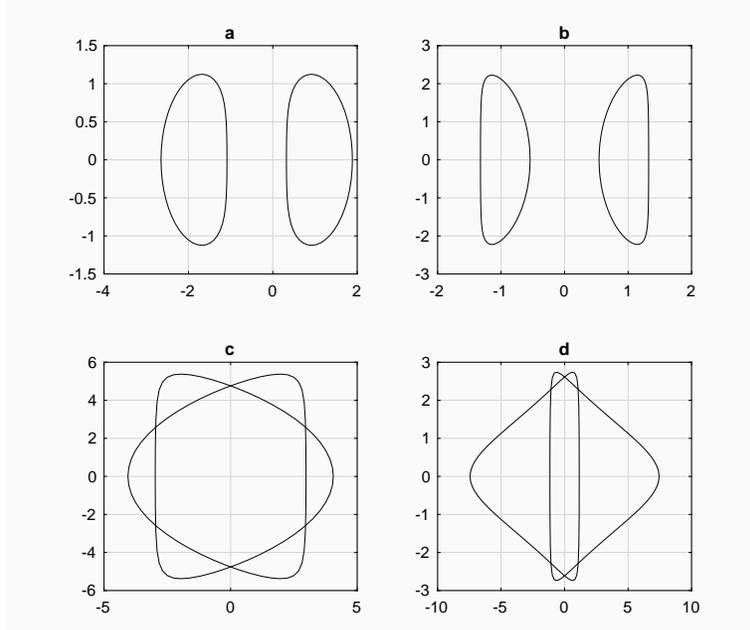}}}
\caption{
Samples of  the contours $L_0$ and $L_1$ when 
$\tau_1/\mu=2$, $\tau_1^\infty/\mu=1 $, $\tau_2=\tau_2^\infty=0$,   $c_{-1}=1$, $r_0=r_1=1$,  $\Gz_0=-1.5$,  $\Gz_1=1.5$, and $\Gz_*=0$.
a: $\Gk_0=\Gk_1=2$, $a_0=\tilde d_0=0$. b-d:  $a_0=-a_1$, and $\tilde d_0=-\tilde d_1$.   b:  $\Gk_0=\Gk_1=0.5$. 
c:  $\Gk_0=\Gk_1=0.9$. b:  $\Gk_0=\Gk_1=1.1$. }
\label{fig2}
\end{figure}

\begin{figure}[t]
\centerline{
\scalebox{0.6}{\includegraphics{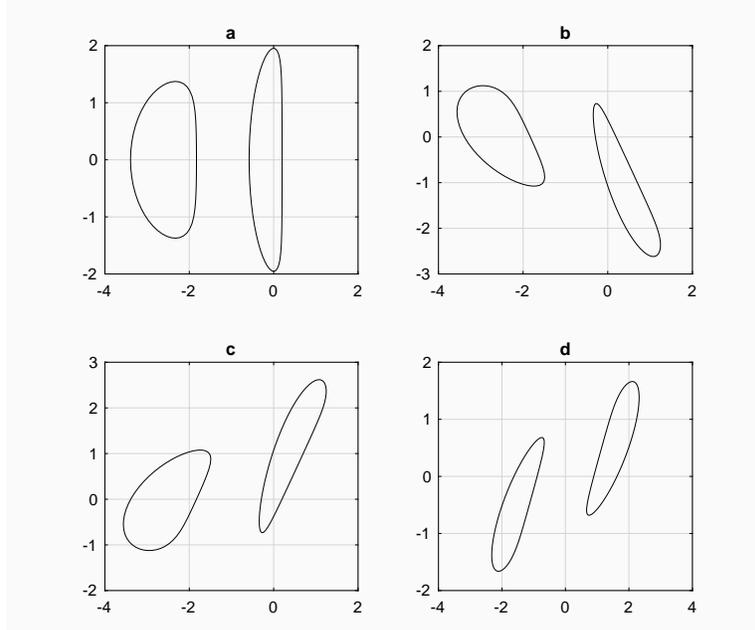}}}
\caption{
Samples of  the contours $L_0$ and $L_1$ when 
$\tau_1/\mu=2$, $\tau_1^\infty/\mu=1 $,  $c_{-1}=1$, $r_0=r_1=1$,  $\Gz_0=-1.5$,  $\Gz_1=1.5$,  $\Gz_*=0$, $a_0=-a_1$, and $\tilde d_0=-\tilde d_1$.
a:  $\tau_2=\tau_2^\infty=0$,  $\Gk_0=2$, $\Gk_1=0.5$.   b:   $\tau_2/\mu=\tau_2^\infty/\mu=1$,  $\Gk_0=2$, $\Gk_1=0.5$. 
c:   $\tau_2/\mu=\tau_2^\infty/\mu=-1$,  $\Gk_0=2$, $\Gk_1=0.5$. d:   $\tau_2/\mu=\tau_2^\infty/\mu=-1$,  $\Gk_0=\Gk_1=10$. 
 }
\label{fig3}
\end{figure}

\begin{figure}[t]
\centerline{
\scalebox{0.6}{\includegraphics{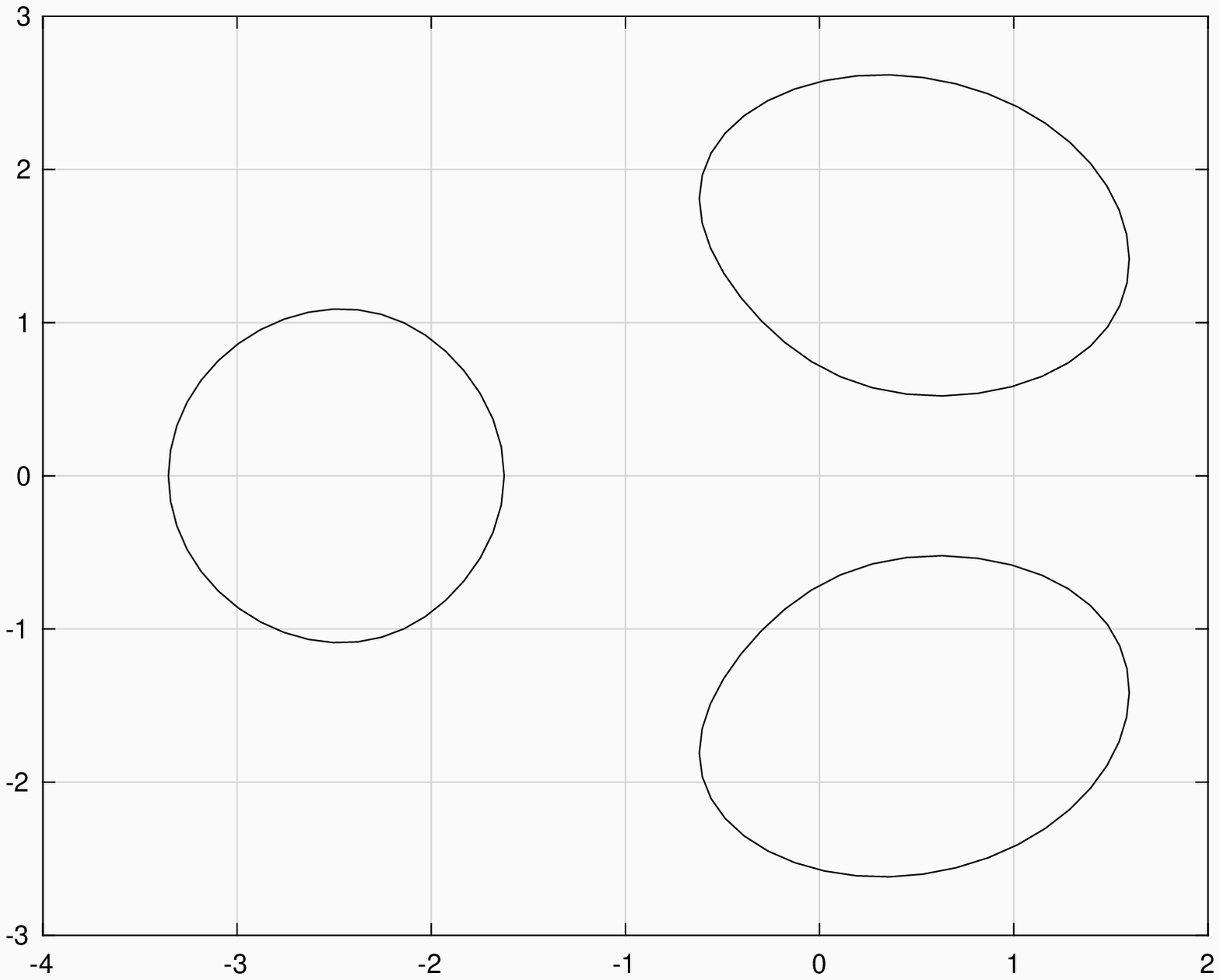}}}
\caption{
Three uniformly stressed inclusions when 
$\tau_1/\mu=2$, $\tau_1^\infty/\mu=1 $,  $\tau_2=\tau_2^\infty=0$, $\Gk_0=\Gk_1=\Gk_2=2$, 
$c_{-1}=1$, $r_0=r_1=r_2=1$,  $\Gz_0=-2$,  $\Gz_1=2e^{\pi i/3}$,  $\Gz_2=2e^{-\pi i/3}$, $\Gz_*=0$, and $a_0=\tilde d_0=0$.
 }
\label{fig4}
\end{figure} 

\begin{figure}[t]
\centerline{
\scalebox{0.6}{\includegraphics{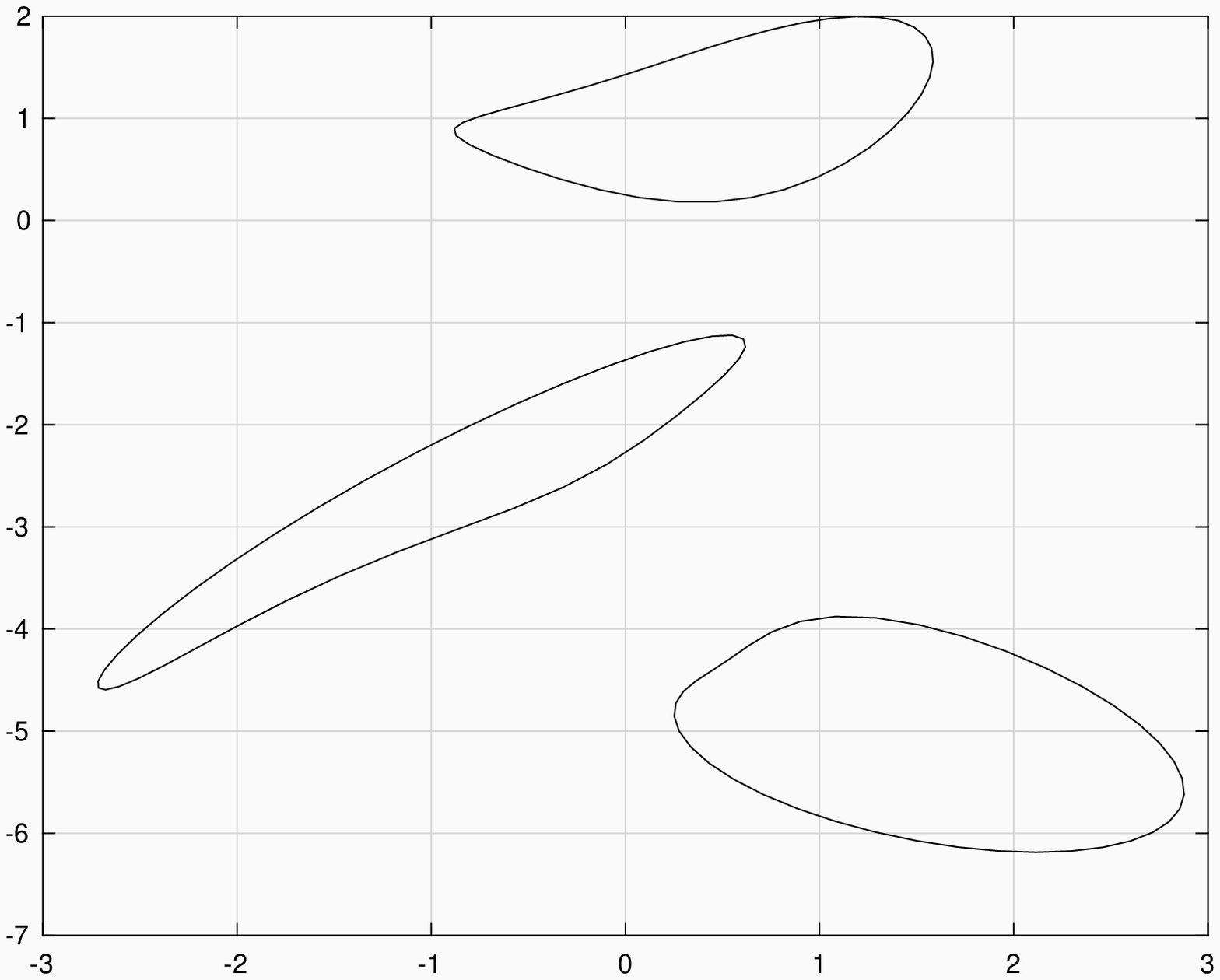}}}
\caption{
Three uniformly stressed inclusions when $\tau_1/\mu=2$, $\tau_1^\infty/\mu=1 $,  $\tau_2/\mu=1$, $\tau_2^\infty/\mu=-1$,
$\Gk_0=2$, $\Gk_1=3$, $\Gk_2=0.5$,
$c_{-1}=1$, $r_0=r_1=r_2=1$,  $\Gz_0=-2$,  $\Gz_1=2e^{\pi i/3}$,  $\Gz_2=2e^{-\pi i/3}$, $\Gz_*=0$, and $a_0=\tilde d_0=0$.
 }
\label{fig5}
\end{figure}

The conformal mapping $\Go(\Gz)$ can be expressed in terms of the solution of the second Riemann-Hilbert problem. From (\ref{3.5}) and (\ref{3.46}) we have
\beq
\Go(\Gx)=c_{-1}\Gx-\fr{\mu}{\bar\tau}\left[g_{2j}^\circ(\Gx)-\tilde d_j+i\Psi_2(\Gx)+i\ov{\Psi_2(T_0(\Gx))}\right], \quad \Gx\in\frak L_j, \quad j=0,1,\ldots, n-1.
\label{4.1}
\eeq
When a point $\Gx$ traverses the circle $\frak L_j$, the point $z=\Go(\Gx)$ describes the circumference of the inclusion $D_j$ ($j=0,1,\ldots,n-1$). 
Formula (\ref{4.1}) found represents  ($3n-4$)-parametric family of conformal mappings of the $n$-connected circular domain into
the $n$-connected physical domain. The free parameters of this family are the radii of the circles $\frak L_j$, $r_j$, ($j=1,2,\ldots,n-1$),
the real center  $\Gz_1$ of the circle $\frak L_1$, and the complex centers $\Gz_j=\Gz_j'+i\Gz_j''$ of the circles
$\frak L_j$ ($j=2,3,\ldots,n-1$). 
Up to transformations of translation, rotation, 
and scaling
the conformal mapping is invariant with respect to the real constants $a_0$ and $\tilde d_0$, the complex parameter $c_{-1}$, and the point $\Gz_*\in \frak D^e$. 

To implement the method, one needs to compute the integrals $\Psi_2(\Gx)$ and $\ov{\Psi_2(T_0(\Gx))}$.  The former integral is
\beq
\Psi_2(\Gx)=\fr{1}{2\pi}\sum_{j=0}^{n-1}\int_{\frak L_j}[g_{2j}^\circ(\Gn)-\tilde d_j]\CK(\Gx,\Gn)d\Gn, \quad \Gx\in \frak L_j, \quad j=0,1,\ldots,n-1.
\label{4.2}
\eeq
The kernel $\CK(\Gx,\Gn)$ can be written in the form
$$
\CK(\Gx,\Gn)=\fr{1}{\Gx-\Gn}-\fr{1}{\Gx-\Gz_*}+\sum_{k_1=0}^{n-1}\sum_{k_2=0, k_2\ne k_1}^{n-1}\left(\fr{1}{\Gn-T_{k_1}T_{k_2}(\Gx)}-\fr{1}{\Gn-T_{k_1}T_{k_2}(\Gz_*)}
\right)
$$ 
$$
+\sum_{k_1=0}^{n-1}\sum_{k_2=0, k_2\ne k_1}^{n-1}\sum_{k_3=0, k_3\ne k_2}^{n-1}\sum_{k_4=0, k_4\ne k_3}^{n-1}
\left(\fr{1}{\Gn-T_{k_1}T_{k_2}T_{k_3}T_{k_4}(\Gx)}
\right.
$$
\beq
\left.
-\fr{1}{\Gn-T_{k_1}T_{k_2}T_{k_3}T_{k_4}(\Gz_*)} \right)+\ldots,\quad \Gx\in \frak L_j, \quad \Gn\in\frak L_l.
\label{4.3}
\eeq
Clearly,  $\CK(\Gx,\Gn)$ is a regular kernel if $j\ne l$, and a singular kernel otherwise. The singular part of the integral (\ref{4.2}) is evaluated
numerically by the formula
$$
\fr{1}{2\pi}\int_{\frak L_l}\fr{\Gf(\Gn)d\Gn}{\Gn-\Gx}=\fr{i}{2(2N+1)}\sum_{j=-N}^N\Gf(\Gz_l+r_le^{i\Gt_j})
$$
\beq
\times\left[1+\fr{2i\sin\fr{N}{2}(\Gt-\Gt_j) \,\sin\fr{N+1}{2}(\Gt-\Gt_j)}{\sin\fr{\Gt-\Gt_j}{2}}\right],\quad \Gx\in\frak L_l,
\label{4.4}
\eeq
where 
\beq
\Gt_j=\fr{2\pi j}{2N+1}, \quad \Gt=-i\ln\fr{\Gx-\Gz_l}{r_l},
\label{4.5}
\eeq  
and $N$ is a sufficiently large positive integer.

For computing the function $\Psi_2(T_0(\Gx))$, we need the kernel $\CK(T_0(\Gx),\Gn)$. Since $T_0(\Gx)=\Gx$, $\Gx\in\frak L_0$, in the case $\Gx\in\frak L_0$
the kernel $\CK(T_0(\Gx),\Gn)$ is given by formula (\ref{4.3}). If $\Gx\in\frak L_j$, $j=1,2,\ldots,n-1$, and $j=l$, then since $\Gs_j T_0(\Gx)=T_j(\Gx)=\Gx$,
we have
$$
    \CK(T_0(\Gx),\Gn)=\fr{1}{\Gn-\Gx}-\fr{1}{\Gn-\Gs_j(\Gz_*)}+\fr{1}{\Gn-T_0(\Gx)}-\fr{1}{\Gn-\Gz_*}
    $$
    \beq
    +\sum_{\Gs\in\frak G\setminus\Gs_0\setminus\Gs_j}
    \left(\fr{1}{\Gn-\Gs T_0(\Gx)}-\fr{1}{\Gn-\Gs(\Gz_*)}\right), \quad \Gx\in\frak L_j.
    \label{4.6}
    \eeq

The inverse model problem of antiplane elasticity solved has $n+4$ parameters, $\tau_1/\mu$,  $\tau_2/\mu$, $\tau_1^\infty/\mu$  $\tau_2^\infty/\mu$, and $\Gk_j$ ($j=0,1,\ldots,n-1$),
The conformal mapping $z=\Go(\Gz)$ possesses $3n+6$ real parameters, $\Gz_j=\Gz_j'+i\Gz_j''$, $r_j$ ($j=0,1,\ldots,n-1$), $c_{-1}=c_{-1}'+ic_{-1}''$, $\Gz_*=\Gz_*'+i\Gz_*''$,
$a_0$, and $\tilde d_0$. Without loss of generality the  ten real parameters $\Gz_0'+i\Gz_0''$, $\Gz_1''$, $r_0$, $c_{-1}'+ic_{-1}''$, $\Gz_*'+i\Gz_*''$,  $a_0$, and $\tilde d_0$  may be arbitrarily
fixed. The other $3n-4$ parameters are free and generate a family of uniformly stressed inclusions.  

Figures \ref{fig1} to \ref{fig3} provides samples of two 
inclusions when the conformal map has two free parameters, the $x$-coordinate of the center and the radius of the circle $\frak L_1$. 
In Fig. \ref{fig1}, we show the profiles of two symmetric inclusions in the case when the loading parameters are
$\tau_1/\mu=2$, $\tau_1^\infty/\mu=1 $, $\tau_2=\tau_2^\infty=0$, the materials of the inclusions and the matrix are characterized by
the parameters $\Gk_1=2$ and $\Gk_2=2$. The circles $\frak L_1$ and $\frak L_2$ are taken to have unit radius, $r_0=r_1=1$, and centered at
$\Gz_0=-1.5$ and $\Gz_1=1.5$, respectively. The point $\Gz_*\in\frak D^e$ is chosen as the origin, $\Gz_*=0$.
It turns out that when the parameters $\Gk_1$ and $\Gk_2$ tend to either 0 or infinity the contours $L_0$ and $L_1$ become slim and, in the limit, they become
two segments. On the other hand, for $\Gk_0$
and $\Gk_1$ close to 1, the contours $L_0$ and $L_1$ may  intersect each other as shown in Fig. \ref{fig2}. Fig. \ref{fig3} provides some examples 
of the  problem and conformal map parameters  which generate nonsymmetric contours $L_0$ and $L_1$.

In the case  of three inclusions, the conformal map has five free parameters. Fig. \ref{fig4} illustrates a symmetric case when $\tau_1/\mu=2$, $\tau_1^\infty/\mu=1 $, 
$\tau_2=\tau_2^\infty=0$, the parameters $\Gk_j$  and the radii $r_j$ are the same: $\Gk_0=\Gk_1=\Gk_2=2$,   $r_0=r_1=r_2=1$, and the circles centers
are taken as $\Gz_0=-2$,  $\Gz_1=2e^{\pi i/3}$,  $\Gz_2=2e^{-\pi i/3}$. We give an example of nonsymmetric  uniformly stressed inclusions in 
 Fig. \ref{fig5} when $\Gk_0=2$, $\Gk_1=3$, $\Gk_2=0.5$, $\tau_1/\mu=2$, $\tau_1^\infty/\mu=1 $,  $\tau_2/\mu=1$, $\tau_2^\infty/\mu=-1$. The radii
 and the centers of the circles $\frak L_j$ are the same as in Fig.  \ref{fig4}.

\setcounter{equation}{0}

\section{Conclusion}

To solve the inverse problem of antiplane elasticity on recovering the shape of $n$ uniformly stressed inclusions, we
proposed to apply the method of conformal mappings from  an $n$-connected external circular domain to the exterior of $n$ inclusions.
The reconstruction of the conformal map requires solving two Schwarz problems on $n$ circles when the right hand-side of the
boundary condition of the second problem is expressed through the solution  of the first Schwarz problem. 
To solve these Schwarz problems, we applied the method of symmetry and linear rational transformations.
This approach brought us to two Riemann-Hilbert problems of the theory of
automorphic functions generated by a Schottky group of symmetric linear rational transformations.
By employing a quasiautomorphic analogue of the Cauchy kernel we derived a series representation
of a family of conformal mappings meeting the requirements of the inverse problem of antiplane elasticity. 
This family possesses $3n-4$ free real parameters which have to be chosen in a sensible manner to avoid 
possible overlapping of the inclusions. 

Numerical results obtained for two- and three-connected domains revealed that when the parameters
$\Gk_j\to 0$ or $\Gk_j\to \infty$, the inclusions tend to transform into segments. Here, $\Gk_j=\mu_j/\mu$, $\mu_j$ and $\mu$ are the shear 
moduli of the inclusions $D_j$ and the matrix $D^e$, and $j=0,1,\ldots,n-1$. 

In this paper the preimage parametric domain was 
chosen to be the exterior of $n$ circles. Another possibility for the parametric domain  is the exterior of $n$ slits
lying in the same line (the real axis for example). In this case the slit map was constructed {\bf(\ref{ant7})} by solving two Riemann-Hilbert
problems on a symmetric genus-$n$ Riemann surface. An advantage of such an approach is its ability to recover
the conformal map by quadratures in the cases of doubly and triply connected domains. However,
when $n\ge 4$, since not each $n$-connected domain $D^e$ can be considered as an image by a slit map
with $n$-slits lying in the same line, the method of slit maps is in general inapplicable. At the same time, if $n\ge 4$ and
the associated Schottky group  of the first class that is the series representation of the quasiautomorphic analogue of the Cauchy kernel 
is absolutely convergent, then the method of circular maps and the Riemann-Hilbert problems of the theory of automorphic functions
works and gives a series representation of the conformal map. The set of domains associated with the first class Schottky group is broader
and includes not only the set of $n$ circles whose centers fall in the same line.

\setcounter{equation}{0}
\section*{Appendix: Single inclusion}\label{A}

Without loss of generality $\frak L_0$ is the unit circle centered at the origin and $a_0=0$.
The solution of the Schwarz problem (\ref{2.10}), (\ref{2.13}) for the unit circle  $\frak L_0$
is given by
$$
F(\Gz)=\Gb_0-i\Gb_1\Gz+i\bar\Gb_1\Gz^{-1},
\eqno{(A.1)}
$$
where
$$
\Gb_1=\Gb_1'+i\Gb_1''=\fr{\bar\tau^\infty-\bar\tau}{\mu}ic_{-1},
\eqno{(A.2)}
$$
and  $c_{-1}$ and $\Gb_0$ are real constants. The solution of the second Schwarz problem (\ref{2.11}), (\ref{2.13}) can be represented in the form
$$
\bar\tau\Go(\Gz)=\Gg_{-1}\Gz^{-1}+\Gg_0+\Gg_1\Gz, 
\eqno{(A.3)}
$$
where $\Gg_j=\Gg_j'+i\Gg_j''$, $j=-1,0,1$. On substituting the expressions (A.3) and (A.1) into
the boundary condition (\ref{2.11}) and replacing $\Gz$ by $e^{i\Gvf}$, $0\le\Gvf\le2\pi$, we derive
$$
\Gg_0'=(\Gb_0-d_0')\fr{\mu_0}{1-\Gk_0}, \quad
\Gg_{-1}'+\Gg_1'=\fr{2\Gb_1''\mu_0}{1-\Gk_0}, \quad
\Gg_{-1}''-\Gg_1''=\fr{2\Gb_1'\mu_0}{1-\Gk_0}.
\eqno{(A.4)}
$$
Finally, by using the second formula in (\ref{2.13})  and the relations (A.4)
we determine the function $\Go(\Gz)$ up to an additive complex constant $\Gg$
$$
\Go(\Gz)=c_{-1}\left(\Gz+\fr{\Gd}{\Gz}\right)+\Gg,
\eqno{(A.5)}
$$
where
$$
\Gd=\fr{2\Gk_0\tau^\infty-(\Gk_0+1)\tau}{(1-\Gk_0)\bar\tau}.
\eqno{(A.6)}
$$
Let $\Gk_0\ne 1$, $\tau\ne 0$, and $\Gd\ne \pm 1$. Then a point $z=\Go(\Gz)$ traces an ellipse $L_0$
whenever the point $\Gz$ traverses the unit circle  $\frak L_0$.

\vspace{.2in}

{\centerline{\Large\bf  References}}

\vspace{.1in}

\begin{enumerate}

\item\label{ria} 
D. Riabuchinsky, Sur la d\'etermination d'une surface d`apr\'es les donn\'ees
qu'elle porte, {\it C.-R. Paris} {\bf 189}  (1929) 629-632.

\item\label{aks1} 
L.A. Aksent'ev, N.B. Il'inskii,  M.T. Nuzhin, 
R.B. Salimov and  G.G. Tumashev,  The theory of inverse boundary value problems for analytic functions and
its applications, {\it Mathematical Analysis, Akad. Nauk SSSR, Vsesoyuz. Inst. Nauchn. i Tekhn. Informatsii  Moscow} {\bf 18} (1980) 67-124.

\item\label{esh} 
J.D. Eshelby,    The determination of the elastic field of an ellipsoidal inclusion, and related problems,
{\it Proc. Roy. Soc. London A} {\bf  241} (1957) 376-396.

\item\label{sen} 
G.P. Sendeckyj,   Elastic inclusion problems in plane elastostatics, {\it Int. J. Solids Structures} {\bf 6} (1970 ) 1535-1543.

\item\label{ru} 
 C.-Q. Ru and P. Schiavone,  On the elliptic inclusion in anti-plane shear,  {\it Math. Mech. Solids} {\bf 1} ( 1996) 327-333.

\item\label{che} 
G.P. Cherepanov,  Inverse problems of the plane theory of elasticity, {\it J. Appl. Math. Mech.} {\bf 38} ( 1974)
915-931. 

\item\label{vig} 
S.B.  Vigdergauz,    Integral equation of the inverse problem of the plane theory of elasticity,
{\it J. Appl. Math. Mech.}  {\bf 40} (1976) 518-522.

\item\label{gra} 
Y. Grabovsky and R.V. Kohn,  Microstructures minimizing the energy of a two phase elastic composite in two space dimensions. II.: The Vigdergauz microstructure, {\it J. Mech. Phys. Solids} {\bf  43} ( 1995 ) 949-972.

\item\label{ant1} 
Y.A. Antipov, Slit maps in the study of equal-strength cavities in n-connected elastic planar domains, {\it SIAM J. Appl. Math.} {\bf 78}  (2018) 320-342.

\item\label{kan} 
H. Kang, E. Kim   and  G.W. Milton,  Inclusion pairs satisfying Eshelby's uniformity property,
{\it SIAM, J. Appl. Math.} {\bf 69} ( 2008)  577-595.

\item\label{wan} 
X. Wang, Uniform fields inside two non-elliptical inclusions, {\it Math. Mech. Solids}  {\bf 17}  (2012) 736-761.

\item\label{liu} 
L.P. Liu, Solutions to the Eshelby Conjectures,  {\it Proc. Roy. Soc. London A} {\bf  464} ( 2008) 573-594.

\item\label{dai} 
M. Dai, C.-Q. Ru and C.-F. Gao,   Uniform strain fields inside multiple
inclusions in an elastic infinite plane under
anti-plane shear,  {\it Math. Mech. Solids}  {\bf 17} (2017) 114-128.

\item\label{ant2}
Y.A. Antipov and V.V. Silvestrov,
Method of Riemann surfaces in the study of supercavitating flow around two hydrofoils in a channel, {\it Physica D} {\bf 235} (2007) 72-81. 

\item\label{ant3}
Y.A. Antipov and V.V. Silvestrov, Circular map for supercavitating flow in a multiply connected domain, {\it Quart. J. Mech. Appl. Math.} {\bf 62} (2009)  167-200.

\item\label{kel}
M.V. Keldysh, 
Conformal mappings of multiply connected domains on canonical domains,
{\it Uspekhi Matem. Nauk} {\bf  6} (1939) 90-119.

\item\label{cou}
R. Courant, {\it Dirichlet's Principle, Conformal Mapping, and Minimal Surfaces} (Interscience Publishers, Inc. New York 1950).

\item\label{for} L.R. Ford, {\it Automorphic Functions}
(McGraw-Hill Book Company, New York 1929).

\item\label{chi} 
L. I. Chibrikova and V. V. Silvestrov, On the question of the effectiveness of the solution of
Riemann's boundary value problem for automorphic functions,  {\it Soviet Math. (Iz. VUZ)} {\bf 12} (1978) 117-121. 
 
\item\label{sil} 
V. V. Silvestrov, The Riemann boundary value problem for symmetric automorphic functions
and its application, {\it Theory of Functions of a Complex Variable and Boundary Value Problems} (Chuvash. Gos. Univ., Cheboksary 1982) 93-107.

\item\label{ant4} Y.A. Antipov and V.V. Silvestrov, Method of automorphic functions in the study of flow around a stack of porous cylinders,  {\it Quart. J. Mech. Appl. Math.} {\bf 60} (2007) 337-366. 

\item\label{ant5} Y.A. Antipov and V.V. Silvestrov, Hilbert problem for a multiply connected circular domain and the analysis of the Hall effect in a plate, 
{\it Quart. Appl. Math.} {\bf 68}  (2010) 563-590.

\item\label{aks2}
 L. A. Aksent'ev, Construction of the Schwarz operator by the symmetry method, {\it Trudy Sem.
Kraev. Zadacham (Kazan)} {\bf 4} (1967) 3-10.

\item\label{ale}
I. A. Aleksandrov and A. S. Sorokin, The problem of Schwarz for multiply connected circular
domains, {\it Siberian Math. J.} {\bf 13} (1973) 671-692.

\item\label{mit} 
V. V. Mityushev and S. V. Rogosin, {\it Constructive Methods for Linear and Nonlinear Boundary
Value Problems for Analytic Functions} (Chapman \& Hall, Boca Raton 2000).

\item\label{kaz}  A. Kazarin and  Y. Obnosov, An exact analytical solution of an $R$-linear conjugation problem for a $n$-phased concentric circular heterogeneous structure,
{\it Appl. Math. Modeling} {\bf 40} (2016) 5292-5300.

\item\label{del} 
T.K. Delillo, A.R. Elcrat and J.A. Pfaltzgraff, Schwarz-Christoffel mappings of multiply connected domains, {\it J. d'Analyse} {\bf 94} (2004) 17-47. 

\item\label{cro} 
D. Crowdy, The Schwarz-Christoffel mapping to multiply-connected polygonal domains, 
{\it Proc. R. Soc.} A {\bf 461} (2005)  2653-2678.

\item\label{ant6} 
Y.A. Antipov and D.G. Crowdy,
Riemann-Hilbert problem for automorphic functions and the Schottky-Klein prime function, {\it Complex Anal. Oper. Theory} 
{\bf 1} (2007) 317-334.

\item\label{bur} 
W. Burnside, On a class of automorphic functions, {\it Proc. London Math. Soc.} {\bf 23} (1892) 49-88.

\item\label{sch} 
F. Schottky, Ueber eine specielle Function, welche bei einer bestimmten linearen Transforma-
tion ihres Arguments unver\"andert bleibt, {\it J. Reine Angew. Math.} {\bf 101} (1887) 227-272.

\item\label{gak}  F. D. Gakhov, {\it Boundary Value Problems} (Pergamon Press, Oxford 1966).

\item\label{ant7} 
Y.A. Antipov, Inverse antiplane problem on $n$ uniformly stressed inclusions, submitted for publication, arXiv:1705.06627v2.

\end{enumerate}

\end{document}